\documentclass[11pt]{article}
\usepackage{latexsym,amssymb}
\usepackage{amsmath, amsbsy}
\usepackage{amsopn, amstext}
\newtheorem{lemma}{Lemma}
\newtheorem{theorem}{Theorem}

\newtheorem{proposition}{Proposition}

\def\NN{\hbox{\rlap{I}\kern.16em N}}
\def\NC{\hbox{\rlap{\kern.24em\raise.1ex\hbox
                  {\vrule height1.3ex width.9pt}}C}}

\title{A contribution to the condition number of the total least squares
problem}
\author{Zhongxiao Jia $^1$, Bingyu Li $^{2,}$\footnote{Corresponding author.}  \\ \\
\small {$^1$ Department of Mathematical Sciences, Tsinghua
University,
Beijing 100084, China} \\
\small {$^2$ School of Mathematics and Statistics, Northeast Normal University, } \\
\small { Changchun 130024, China } \\
\small {jiazx@tsinghua.edu.cn (Z. Jia)  \,\,  mathliby@gmail.com (B.
Li)} }

\date{}

\begin{document}

\maketitle

\begin{abstract}
This paper concerns cheaply computable formulas and bounds for the
condition number of the TLS problem. For a TLS problem with data
$A$, $b$, two formulas are derived that are simpler and more compact
than the known results in the literature. One is derived by
exploiting the properties of Kronecker products of matrices.
The other is obtained by making use of the singular value
decomposition (SVD) of $[A \,\,b]$, which allows us to compute the
condition number cheaply and accurately. We present lower and upper
bounds for the condition number that involve the
singular values of $[A \,\, b]$ and the last entries of the right
singular vectors of $[A \,\, b]$. We prove that they are always
sharp and can estimate the condition number accurately by no more
than four times. Furthermore, we establish a few other lower and
upper bounds that involve only a few singular values of $A$ and $[A
\,\, b]$. We discuss how tight the bounds are. These bounds are
particularly useful for large scale TLS problems since for them any
formulas and bounds for the condition number involving all the
singular values of $A$ and/or $[A \ b]$ are too costly to be
computed. Numerical experiments illustrate that our bounds are
sharper than a known approximate condition number in the literature.

\end{abstract}

\vskip 5pt \noindent {\bf Keywords:}  total least squares, condition
number, singular value decomposition.

\vskip 5pt \noindent {\bf AMS subject classification (2000):}
65F35.

\section{Introduction}

 For given $A \in \mathbb{R}^{m \times n} (m
> n)$, $b \in \mathbb{R}^m$, the TLS problem can be formulated as
(see, e.g., \cite{GolubVanLoan:1980,PaigeStrakos:2002})
\begin{equation}\label{ScaledTLSprob}
  \min \|[E\,\,r]\|_F,
  \text{\quad subject to  \quad}
   b - r \in \mathcal{R}(A+E)  ,
\end{equation}
where $\|\cdot\|_F$ denotes the Frobenius norm of a matrix and
$\mathcal{R}(\cdot)$ denotes the range space. Suppose that $[E_{TLS}
\,\,r_{TLS}]$ solves the above problem. Then $x = x_{TLS}$ that
satisfies the equation $(A+E_{TLS}) x = b - r_{TLS}$ is called the
TLS solution of (\ref{ScaledTLSprob}).

The condition number measures the worst-case sensitivity of a
solution of a problem to small perturbations in the input data.
Combined with backward errors, it provides an approximate local
linear upper bound on the computed solution. Since the 1980's,
algebraic perturbation analysis for the TLS problem has been studied
extensively; see
\cite{FierroBunch:1996,GolubVanLoan:1980,Liu:1996,WeiM:1992,WeiM:1998}
and the references therein. In recent years, asymptotic perturbation
analysis and TLS condition numbers have been studied; see, e.g.,
\cite{BaboulinGratton:2010,LiJia:2009,ZhouLinWeiQiao:2008}.

In the present paper, we continue our work in \cite{LiJia:2009} that
studied the condition number of the TLS problem. We will derive a
number of results. Firstly, we establish two formulas that are
simpler and more suitable for computational purpose than the known
results in the literature. One is derived by exploiting the
properties of Kronecker products of matrices. It improves the
formulas given in \cite{LiJia:2009,ZhouLinWeiQiao:2008}, is
independent of Kronecker products of matrices and makes its
computation convenient. The other is obtained by making use of
the SVD of $[A \,\,b]$, which can be used to compute the condition
number more cheaply and accurately than that in
\cite{BaboulinGratton:2010}. Secondly, we present lower and upper
bounds for the condition number that involve the
singular values of $[A \,\, b]$ and the last entries of the right
singular vectors of $[A \,\, b]$. We prove that these bounds are
always sharp and can estimate the condition number accurately by no
more than four times. Finally, we focus on cheaply computable
bounds of the TLS condition number. We establish lower and upper
bounds that involve only a few singular values of $A$ and $[A \,\,
b]$. We discuss how tight the bounds are. These bounds are
particularly useful for large scale TLS problems since for them any
formulas and bounds for the condition number involving all the
singular values of $A$ and/or $[A \ b]$ are too costly to be
computed. So we can compute these bounds by calculating only a few
singular values of $A$ and/or $[A\ b]$ using some iterative solvers
for large SVDs. In \cite{BjorckHeggernesMatstoms:2000}, an
approximate TLS condition number is presented and is applied to
evaluate conditioning of the TLS problem there. In this paper, we
present numerical experiments to demonstrate a possibly great
improvement of one of our upper bounds over the approximate
condition number in \cite{BjorckHeggernesMatstoms:2000}.

The paper is organized as follows. In Section~2, we give some
preliminaries necessary. In Section 3, we present computable
formulas of the TLS condition number. The straightforward bounds on
the TLS condition number are considered in Section 4. In Section 5,
we present numerical experiments to show the tightness of bounds for
the TLS condition number. We end the paper with some concluding
remarks in Section 6.

Throughout the paper, for given positive integers $m, n$,
$\mathbb{R}^n$ denotes the space of $n$-dimensional real column
vectors, $\mathbb{R}^{m \times n}$ denotes the space of all $m
\times n$ real matrices.  $\|\cdot\|$ and $\|\cdot\|_F$ denote
2-norm and Frobenius norm of their arguments, respectively.  Given a
matrix $A$, $A(1:i,1:i)$ is a Matlab notation that denotes the $i$th
leading principal submatrix of $A$, and $\sigma_i(A)$ denotes the
$i$th largest singular value of $A$. For a vector $a$, $a(i)$
denotes the $i$th component of $a$, and ${\rm{diag}}(a)$ is a
diagonal matrix whose diagonals are given as components of $a$.
$I_n$ denotes the $n \times n$ identity matrix, $O_{mn}$ denotes the
$m \times n$ zero matrix, whereas $O$ denotes a zero matrix whose
order is clear from the context. For matrices $A =[a_1, \ldots, a_n]
= [A_{ij}] \in \mathbb{R}^{m \times n} $ and $B$, $A \otimes B
=[A_{ij} B] $ is the Kronecker product of $A$ and $B$, the linear
operator ${\rm{vec}}: \mathbb{R}^{m \times n} \rightarrow
\mathbb{R}^{mn} $ is defined by ${\rm{vec}}(A) = [a^T_1, \ldots,
a^T_n]^T$ for $A \in \mathbb{R}^{m \times n}$.

\section{Preliminaries}

Throughout the paper,  we let $\hat{U}^T A \hat{V} =
{\rm{diag}}(\hat{\sigma}_1, \ldots, \hat{\sigma}_n)$ be the thin SVD
of $A \in \mathbb{R}^{m \times n}$, where  $\hat{\sigma}_1 \geq
\cdots \geq \hat{\sigma}_n$, $\hat{U}\in \mathbb{R}^{m \times n}$,
$\hat{U}^T \hat{U} = I_{n}$, $\hat{V} \in \mathbb{R}^{n \times n}$,
$\hat{V}^T \hat{V} = I_{n}$. Let $U^T [A \,\, b] V =
{\rm{diag}}(\sigma_1, \ldots, \sigma_{n+1})$ be the thin SVD of $[A
\,\, b] \in \mathbb{R}^{m \times (n+1)}$, where $\sigma_1 \geq
\cdots \geq \sigma_{n+1}$, $U = [u_1, \ldots, u_{n+1}]\in
\mathbb{R}^{m \times (n+1)}$,  $U^T U = I_{n+1}$, $V = [v_1, \ldots,
v_{n+1}] \in \mathbb{R}^{(n+1) \times (n+1)}$, $V^T V = I_{n+1}$.

The following result presents an existence and uniqueness condition
for the TLS solution \cite{GolubVanLoan:1980}.
\begin{theorem}\label{Golub-vancondition}
 If \begin{equation}\label{Golub-vancondition-inequality}
\sigma_{n+1} < \hat{\sigma}_{n},
\end{equation}then the TLS problem (\ref{ScaledTLSprob}) has a
unique solution $x_{TLS}$. Moreover,
\begin{eqnarray}\label{STLSsolutionExp}
 x_{TLS} & = & (A^T A - \sigma^2_{n+1} I )^{-1} A^T b \\ \label{SVDSTLSsol}
 &=& -
\left[\frac{v_{n+1}(1)}{v_{n+1}(n+1)}, \ldots,
\frac{v_{n+1}(n)}{v_{n+1}(n+1)} \right]^T.
\end{eqnarray}
\end{theorem}

In the paper, it is always assumed that condition
(\ref{Golub-vancondition-inequality}) holds. We note that, for a
given TLS problem (\ref{ScaledTLSprob}), if $\sigma_{n+1} =0$, then
$b\in \mathcal{R}(A)$. In this case, the system of equations $A x =
b$ is compatible, and we can take $[E \,\, r ] = O$. As in
\cite{LiJia:2009,ZhouLinWeiQiao:2008}, in the sequel, we do not
consider this trivial case and assume that
\begin{equation}\label{assumption}
0 < \sigma_{n+1} < \hat{\sigma}_{n}.
\end{equation}

We will use the following properties of the TLS problem, which are
in \cite{GolubVanLoan:1980}:
\begin{equation}\label{optimalfunction}
\sigma^2_{n+1} = \frac{\|r\|^2}{1 + \|x\|^2}
\end{equation}and
\begin{equation}\label{gradienteqn}
A^T r = \frac{\|r\|^2}{1 + \|x\|^2} x = \sigma^2_{n+1} x,
\end{equation}where $x=x_{TLS}$, $r = A x - b$. By (\ref{SVDSTLSsol}), it
holds that
\begin{equation}\label{solusv}
v_{n+1} = \frac{1}{ \sqrt{1 +\|x\|^2}} \left[
                                                 \begin{array}{c}
                                                    x \\
                                                   -1 \\
                                                 \end{array}
                                               \right]
\end{equation}
up to a sign $\pm 1$.

The following basic properties of the Kronecker products of matrices
are needed later and can be found in \cite{Graham:1981}:
\begin{eqnarray}\nonumber
&& (A \otimes C)(B\otimes D) = (A B) \otimes (C D), \\ \nonumber &&
(A \otimes B)^T = A^T \otimes B^T,
\end{eqnarray}where $A, B, C, D$ are matrices with appropriate
sizes.

\section{Computable formulas for the TLS condition number}

Throughout the paper, we follow the definition of condition number
in \cite{GohbergKoltracht:1993,Rice:1966}. Let $g: \mathbb{R}^p
\longrightarrow \mathbb{R}^q$ be a continuous map in normed linear
spaces defined on an open set $D_g \subset \mathbb{R}^p$. For a
given $a \in D_g$, $a \neq 0$, with $g(a) \neq 0$, if $g$ is
differentiable at $a$, then the relative condition number of $g$ at
$a$ is
\begin{equation}\label{Rice_rel}
\kappa^r_g(a) =  \frac{\|g'(a)\| \|a\|}{\|g(a)\|}
\end{equation}
and the absolute condition number of $g$ is
\begin{equation}\label{Rice_abs}\kappa_g(a) = \|g'(a)\|,
\end{equation}
where $g'(a)$ denotes the derivative of $g$ at $a$. Given the TLS
problem (\ref{ScaledTLSprob}), let $\tilde{A}= A + \Delta A$,
$\tilde{b} = b+ \Delta b$, where $\Delta A$ and $\Delta b$ denote
the perturbations in $A$ and $b$, respectively. Consider the
perturbed TLS problem
\begin{equation}\label{perturbedScaledTLS}
\min \|[E \,\,r]\|_F \,\,\,{\text{subject to}} \,\, \tilde{b}-r \in
\mathcal{R}(\tilde{A} + E).
\end{equation}

In \cite{LiJia:2009}, we have established the following result.

\begin{theorem}\label{MyProp}
Suppose the TLS problem (\ref{ScaledTLSprob}) satisfies
(\ref{assumption}). Denote  by $x = x_{TLS}$ the TLS solution, and
define $r = A x - b$, $G(x) = [x^T \,\,-1] \otimes I_m$. If
$\|[\Delta A \,\, \Delta b]\|_F$ is small enough, then the perturbed
problem (\ref{perturbedScaledTLS}) has a unique TLS solution
$\tilde{x}_{TLS}$. Moreover,
\begin{equation}\label{vecExpression}
\tilde{x}_{TLS} = x_{TLS} + K~ \left[
                                                       \begin{array}{c}
                                                         {\rm{vec}}(\Delta A) \\
                                                         \Delta b \\
                                                       \end{array}
                                                     \right]
+ \mathcal{O}(\|[\Delta A \,\, \Delta b]\|^2_F),
\end{equation}
where \begin{equation}\label{bigK} K = \left(A^T A - \sigma^2_{n+1}
I_{n} \right)^{-1} \left( 2 A^T \frac{r}{\|r\|} \frac{r^T}{\|r\|}
G(x) - A^T G(x) - [I_n \otimes r^T \,\,O] \right).
\end{equation}
\end{theorem}

Denote $a ={\rm{vec}}(A)$. Based on Theorem \ref{MyProp}, in a small
neighborhood of $[a^T,b^T]^T \in \mathbb{R}^{m (n+1)}$,  we can
define the function
\begin{equation}\nonumber
\begin{array}{ccc}
  g: \mathbb{R}^{m (n+1)}  &  \longrightarrow & \mathbb{R}^n \\
  \small{\left[
    \begin{array}{c}
      \tilde{a} \\
      \tilde{b} \\
    \end{array}
  \right]}
   & \longmapsto & \tilde{x}= (\tilde{A}^T \tilde{A} -
   \tilde{\sigma}^2_{n+1}I_n)^{-1} \tilde{A}^T
  \tilde{b},
\end{array}
\end{equation}
where $\tilde{a} = a + {\rm{vec}}(\Delta A)={\rm{vec}}( \tilde{A})
$, $\tilde{b} = b+ \Delta b$, and $\tilde{x}$ is the solution to the
perturbed TLS problem (\ref{perturbedScaledTLS}). In particular,
$g([a^T, b^T]^T) = x$. Thus, we have the following theorem.

\begin{theorem}\label{knoneckercondform}
Let $\kappa_{g}(A,b)$ and $\kappa^r_{g}(A,b)$ be the absolute and
relative condition numbers of the TLS problem, respectively. Then
\begin{equation}\label{knoneckercond}
\kappa_{g}(A,b) = \|K\|, \,\,\kappa^r_{g}(A,b) = \frac{\|K\| \|[A
\,\,b]\|_F}{\|x\|},
\end{equation} where $K$ is defined as in (\ref{bigK}).
\end{theorem}
{\bf{{Proof}}}.\,\, By the definition of $g$ and
(\ref{vecExpression}), we see that $g$ is differentiable at
$[a^T,b^T]^T$ and $g'\left([a^T,b^T]^T \right) = K$. Then the
assertion follows from (\ref{Rice_rel}) and (\ref{Rice_abs}). \hfill
$\Box$

The dependence of Kronecker product of matrices for $K$ makes the
computation of $\kappa_{g}(A,b)$ via (\ref{knoneckercond}) too
costly. The same are the formulas given in
\cite{LiJia:2009,ZhouLinWeiQiao:2008}. For a computational purpose,
we will present a new formula of $\kappa_{g}(A,b)$ that has a
simpler and clearer form. To this end, we need a lemma.

\begin{lemma}\label{forCholeskyfac}
Let \begin{equation}\nonumber C = A^T A + \sigma^2_{n+1} I_n -
\frac{2 \sigma^2_{n+1} x x^T}{\|x\|^2+ 1}.
\end{equation} Then $C$ is positive
definite.
\end{lemma}
{\bf{{Proof}}}.\,\, Noticing that
\begin{equation}\label{mytransform}
 C =   A^T A - \sigma^2_{n+1} I_n
+2 \sigma^2_{n+1} \left(I_n - \frac{x x^T}{1 + \|x\|^2} \right),
\end{equation} and that $A^T A - \sigma^2_{n+1} I_n$ and $I_n - \frac{x x^T}{1 + \|x\|^2}$ are both positive definite,
 we complete the proof of the lemma. \hfill $\Box$

\begin{theorem}\label{treatkronecker}
Let  $ A^T A + \sigma^2_{n+1} I_n - \frac{2 \sigma^2_{n+1} x
x^T}{\|x\|^2+ 1} = L L^T$ be the Cholesky factorization. Then
\begin{equation}\label{choleskyformcond}
\kappa_{g}(A,b) = \sqrt{\|x\|^2+ 1} \left\| (A^T A -
\sigma^2_{n+1}I_n)^{-1} L \right\|.
\end{equation}
\end{theorem}
{\bf{{Proof}}}. Consider expression (\ref{bigK}) of $K$. By the
properties of Kronecker product of matrices, we get

\begin{equation}\nonumber
G(x)  G^T(x) = \left( [ x^T \,-1] \otimes I_m  \right)  \left(
\left[
                                \begin{array}{c}
                                  x \\
                                  -1 \\
                                \end{array}
                              \right] \otimes I_m  \right)
= (  \|x\|^2 +  1 ) I_m,
\end{equation}
\begin{equation}\nonumber
[I_n \otimes r^T \,O]  G^T(x) = [ I_n \otimes r^T \,O] \left[
                            \begin{array}{c}
                               x \otimes I_m \\
                              - I_m \\
                            \end{array}
                          \right]
= (I_n \otimes r^T) ( x \otimes I_m) =
 x r^T
\end{equation}
and
\begin{equation}\nonumber
[I_n \otimes r^T \,O]  \left[
                                      \begin{array}{c}
                                        I_n \otimes r \\
                                        O \\
                                      \end{array}
                                    \right]
 =(I_n \otimes r^T) (I_n \otimes r) = \|r\|^2 I_n.
\end{equation}
Thus, we have
\begin{eqnarray}\nonumber
&& \left( 2 A^T \frac{r}{\|r\|} \frac{r^T}{\|r\|}G(x)  - A^T G(x) -
[I_n \otimes r^T \,O]  \right) \\ \nonumber && \cdot \left( 2
G^T(x)\frac{r}{\|r\|} \frac{r^T}{\|r\|} A -  G^T(x) A -
 \left[\begin{array}{c}
 I_n \otimes r \\
 O \\
 \end{array}
 \right] \right) \\\nonumber
&=& (  \|x\|^2+ 1  ) A^T A + \|r\|^2 I_n - x r^T A - A^T r x^T \\
\nonumber &=&  ( \|x\|^2+ 1  ) A^T A + \|r\|^2 I_n- 2 \sigma^2_{n+1}
x x^T.
 \end{eqnarray}
The last equality used $A^T r x^{T}= \sigma^2_{n+1} x x^T$, which
follows from (\ref{gradienteqn}). Denote $P = A^T A - \sigma^2_{n+1}
I_n$. We get
\begin{eqnarray}\label{firstchange}
 K  K^T &=& P^{-1} \left( (  \|x\|^2+ 1 ) A^T A +
\|r\|^2 I_n- 2 \sigma^2_{n+1} x x^T   \right) P^{-1}
\\ \label{ForWK}
&=& (\|x\|^2+ 1) P^{-1} \left(  A^T A + \sigma^2_{n+1} I_n - \frac{2
\sigma^2_{n+1} x x^T}{\|x\|^2+ 1} \right) P^{-1}.
\end{eqnarray}In the last equality we used (\ref{optimalfunction}).
From Theorem \ref{knoneckercondform}, we have
\begin{equation}\nonumber
\kappa_{g}(A,b) =(\|x\|^2+ 1)^{\frac{1}{2}} \left\|P^{-1} \left(
A^T A + \sigma^2_{n+1} I_n - \frac{2 \sigma^2_{n+1} x x^T}{\|x\|^2+
1} \right) P^{-1} \right\|^{\frac{1}{2}}. \label{similartoBaboulin}
\end{equation}
Based on Lemma \ref{forCholeskyfac}, we complete the proof. \hfill
$\Box$

Compared with the formula of $\kappa_{g}(A,b)$ in Theorem
\ref{knoneckercondform}, the formula in Theorem \ref{treatkronecker}
does not involve the Kronecker product of matrices and makes its
computation convenient. However, if $\hat{\sigma}_n$ and
$\sigma_{n+1}$ are close, then $A^T A - \sigma^2_{n+1}I_n$ becomes
ill conditioned. Therefore, it may be hard to use
(\ref{choleskyformcond}) to calculate $\kappa_{g}(A,b)$ accurately.
Next we derive a new formula that can be used to compute the
condition number accurately.

\begin{theorem}\label{myformationforcond}
Let $U^T[A \,\,b] V = {\rm{diag}}(\sigma_1, \ldots, \sigma_{n+1})$
be the SVD of $[A \,\,b]$ with $V = [v_1, \ldots, v_{n+1}]$.
Denote $V_{11} = V(1:n, 1:n)$. Then
\begin{equation}\label{mynewclosed}
\kappa_{g}(A,b) =\sqrt{ \|x\|^2 + 1} ~\|V^{-T}_{11} S\|,
\end{equation}
where $S = {\rm{diag}}([s_1, \ldots, s_n])$, $s_i =
\frac{\sqrt{\sigma^2_i + \sigma^2_{n+1}}}{\sigma^2_i -
\sigma^2_{n+1}}$, $i = 1, \ldots, n$.
\end{theorem}

{\bf{Proof.}} Denote $P = A^T A - \sigma^2_{n+1}  I_n$. From
(\ref{ForWK}),  we have
\begin{equation}\label{forcond}
 \frac{1}{\|x\|^2 + 1 }~K K^T = P^{-1} + 2 \sigma^2_{n+1} P^{-1} \left( I_n  -
 \frac{ x x^T}{1 + \|x\|^2} \right) P^{-1}.
\end{equation}Note that
\begin{eqnarray}\nonumber
[A \,\,b]^T [A \,\,b] - \sigma^2_{n+1} I_{n+1}& =& \sum^{n+1}_{i=1}
\sigma^2_i v_i v^T_i - \sigma^2_{n+1} \sum^{n+1}_{i=1} v_i v^T_i \\
\nonumber &=& \sum^{n}_{i=1} (\sigma^2_i -\sigma^2_{n+1} )v_i v^T_i.
\end{eqnarray}
We get
\begin{eqnarray}\nonumber
P = A^T A - \sigma^2_{n+1} I_{n} &=& [I_n \,\,0] \sum^{n}_{i=1}
(\sigma^2_i -\sigma^2_{n+1} )v_i v^T_i \left[
\begin{array}{c}
I_n \\
0 \\
\end{array}
\right] \\ \nonumber &=& [I_n \,\,0] [v_1, \ldots, v_n] \left[
                                     \begin{array}{ccc}
                                       \sigma^2_1 - \sigma^2_{n+1} &  &  \\
                                        & \ddots &  \\
                                        &  & \sigma^2_n - \sigma^2_{n+1} \\
                                     \end{array}
                                   \right] \left[
                                             \begin{array}{c}
                                               v^T_1 \\
                                               \vdots \\
                                               v^T_n \\
                                             \end{array}
                                           \right]  \left[
                                                      \begin{array}{c}
                                                        I_n \\
                                                        0 \\
                                                      \end{array}
                                                    \right] \\ \label{forP}
&=& V_{11} \left[
                                     \begin{array}{ccc}
                                       \sigma^2_1 - \sigma^2_{n+1} &  &  \\
                                        & \ddots &  \\
                                        &  & \sigma^2_n - \sigma^2_{n+1} \\
                                     \end{array}
                                   \right] V^T_{11} := V_{11} \Lambda V^T_{11}.
                                   \label{pdecom}
\end{eqnarray}
Similarly,  by (\ref{solusv}), since $v_{n+1} = \frac{1}{\sqrt{1 +
\|x\|^2}} \left[
                                          \begin{array}{c}
                                            x \\
                                            -1 \\
                                          \end{array}
                                        \right],
$ we have
\begin{eqnarray}\nonumber
 I_{n+1} - \frac{1}{{1 + \|x\|^2}} \left[
 \begin{array}{cc}
 x x^T & -x \\
 -x & 1 \\
 \end{array}
 \right] = I_{n+1} - v_{n+1} v^T_{n+1} =[v_1, \ldots, v_n] \left[
 \begin{array}{c}
 v^T_1 \\
 \vdots \\
 v^T_n \\
 \end{array}
 \right]
\end{eqnarray}and
\begin{equation}\label{forV11}
 I_n - \frac{x x^T}{1 + \|x\|^2 } = V_{11} V^T_{11}.
\end{equation}
By (\ref{forV11}), we see that $V_{11}$ is invertible. Combining
(\ref{forP}) and (\ref{forV11}), we have
\begin{eqnarray}\nonumber
&& P^{-1} + 2 \sigma^2_{n+1} P^{-1} \left( I_n  -  \frac{ x x^T}{1 +
\|x\|^2} \right) P^{-1} \\ \nonumber &=& V^{-T}_{11} \Lambda^{-1}
V^{-1}_{11} + 2 \sigma^2_{n+1} \left( V^{-T}_{11} \Lambda^{-1}
V^{-1}_{11} \right)  V_{11} V^T_{11} \left( V^{-T}_{11}
\Lambda^{-1}V^{-1}_{11} \right) \\ \label{Wewilluseit} &=&
V^{-T}_{11} \Lambda^{-1} V^{-1}_{11} + 2 \sigma^2_{n+1} V^{-T}_{11}
\Lambda^{-2} V^{-1}_{11}
\\ \label{proofused} &=& V^{-T}_{11} \left( \Lambda^{-1} + 2
\sigma^2_{n+1} \Lambda^{-2} \right) V^{-1}_{11} = \left( V^{-T}_{11}
S \right) \left( V^{-T}_{11} S \right)^T.
\end{eqnarray}Then by (\ref{forcond}) and Theorem \ref{knoneckercondform}
we get the desired equality.
 \hfill
$\Box$

By Theorem \ref{myformationforcond}, we can calculate
$\kappa_{g}(A,b)$ by solving a linear system with the coefficient
matrix $V^T_{11}$. Next we show what the condition number of
$V^T_{11}$ is exactly.

\begin{theorem}\label{svforV11}
For $V_{11}$, we have
\begin{equation}\label{V11singularvalues}
\sigma_1(V_{11}) = 1, \ldots, \sigma_{n-1}(V_{11}) = 1,
\sigma_n({V_{11}}) = \frac{1}{\sqrt{1 + \|x\|^2}}
\end{equation}
and
\begin{equation}\label{V11conditionnumber}
\kappa(V^T_{11}) = \frac{\sigma_1(V_{11})}{\sigma_n(V_{11})}=
\sqrt{1 + \|x\|^2}.
\end{equation}
\end{theorem}
{\bf{Proof.}} By the definition of $V_{11}$ and the interlacing
property \cite[p.103]{Wilkinson:1965} for eigenvalues of symmetric
matrices, we get
\begin{equation}\nonumber
\sigma_1(V_{11}) = 1, \ldots, \sigma_{n-1}(V_{11}) = 1.
\end{equation}
Noticing that
\begin{equation}\nonumber
V_{11} V^T_{11} x = \left(I_n - \frac{1}{1 + \|x\|^2} x x^T \right)
x = x - \frac{\|x\|^2}{1 + \|x\|^2} x = \frac{1}{1 + \|x\|^2} x,
\end{equation}
we know that $\frac{1}{1 + \|x\|^2}$ is an eigenvalue of $V_{11}
V^T_{11}$, that is, $\sigma_n({V_{11}}) = \frac{1}{\sqrt{1 +
\|x\|^2}}$. Thus, we have proved  (\ref{V11singularvalues}) and
(\ref{V11conditionnumber}). \hfill $\Box$

A different SVD-based closed formula for $\kappa_{g}(A,b)$ appears
in \cite{BaboulinGratton:2010}. It is shown in
\cite{BaboulinGratton:2010} that
\begin{equation}\label{Baboulin}
\kappa_{g}(A,b) =\sqrt{ \|x\|^2 + 1} \left\| \hat{D}~[\hat{V}^T
\,\,O_{n,1}]~  V~ [D \,\,\,O_{n,1}]^T \right\|,
\end{equation}where
\begin{eqnarray}\nonumber
\hat{D} &= & {\rm{diag}}\left( \left[ \frac{1}{\hat{\sigma}^2_1 -
\sigma^2_{n+1} }, \ldots, \frac{1}{\hat{\sigma}^2_n - \sigma^2_{n+1}
} \right] \right), \\ \nonumber D &=& {\rm{diag}} \left( \left[
\sqrt{\sigma^2_1 + \sigma^2_{n+1}}, \ldots, \sqrt{\sigma^2_n +
\sigma^2_{n+1}} \right]\right).
\end{eqnarray}
Compared with (\ref{Baboulin}), our (\ref{mynewclosed}) is simpler
and more compact. Furthermore,  (\ref{Baboulin}) depends on the
singular values and right singular vectors of both $A$ and $[A
\,\,b]$. In contrast, (\ref{mynewclosed}) involves only singular
values and right singular vectors of $[A \,\,b]$. Therefore, the
computational cost of the condition number by (\ref{mynewclosed}) is
half of that by (\ref{Baboulin}). Furthermore, the following
example shows that the computed results by (\ref{mynewclosed}) can
be more accurate than those by (\ref{Baboulin}).

{\bf{A small example.}} We construct a TLS problem with
$\hat{\sigma}_n$ and $\sigma_{n+1}$ very close. We generate $A, b$
by $[A \,\,b] = {\text{\bf{generate}}}Ab\alpha(m,n, \alpha)$ (see
Appendix) by taking $m=15, n=10$, $\alpha = 10^{-8}$.
 \begin{table}[h]\label{tab1}
\begin{center}
{ \doublerulesep18.0pt \tabcolsep0.01in
\begin{tabular}{|c  | c | c |c| c |}
\hline   $~\sigma_{n+1}/\sigma_n~$ &
$~\sigma_{n+1}/\hat{\sigma}_n~$&
 $~{\kappa}^r_g(\ref{knoneckercond})~$  &
$~{\kappa}^r_g(\ref{mynewclosed})~$& $~\kappa^r_g(\ref{Baboulin})~$
 \\
\hline
 $0.608$ &
 $~ 1-1.49\times 10^{-15}$ & $~-~$    &$1.13\times 10^9$ & $3.12\times 10^8$
 \\
   \hline
\end{tabular}}
\end{center}
\end{table}

In the table, $~\sigma_{n+1}/\sigma_n~$ and
$~\sigma_{n+1}/\hat{\sigma}_n~$ denote the quotients of
$\sigma_{n+1}$ over $\sigma_n$ and $\hat{\sigma}_n$, respectively.
${\kappa}^r_g(\ref{knoneckercond})$,
${\kappa}^r_g(\ref{mynewclosed})$ and $\kappa^r_g(\ref{Baboulin})$
denote the computed $\kappa^r_g(A,b)$ by calculating $\kappa_g(A,b)$
via (\ref{knoneckercond}), (\ref{mynewclosed}) and (\ref{Baboulin}),
respectively. $\sigma_{n+1}$ and $\hat{\sigma}_n$ being so close
makes $A^T A - \sigma^2_{n+1} I_n$ numerically singular and makes
${\kappa}^r_g(\ref{knoneckercond})$ unreliable completely, so the
result of ${\kappa}^r_g(\ref{knoneckercond})$ is omitted.

We comment that ${\kappa}^r_g(\ref{mynewclosed})$ is reliable as, by
the remark in Appendix and Theorem \ref{svforV11}, $\kappa(V^T_{11})
= \alpha^{-1} =10^8$. This means that computing $\kappa_g(A,b)$ via
(\ref{mynewclosed}) amounts to solving a moderately ill-conditioned
linear system. Furthermore, the right-hand side $S$ of the system
can be constructed with high accuracy since $\sigma_{n+1}$ and
$\sigma_i$ are not close: $\frac{\sigma_{n+1}}{\sigma_i} \leq
\frac{\sigma_{n+1}}{\sigma_n} = 0.608$, $i=1,2,\ldots, n$. In
contrast, $\kappa^r_g(\ref{Baboulin})$ is inaccurate since computing
$\kappa_g(A,b)$ via (\ref{Baboulin}) involves the diagonal matrix
$\hat{D}$ and the closeness of $\sigma_{n+1}$ (about $0.299$) and
$\hat{\sigma}_n$ makes its last diagonal entry both very large
(about $10^{15}$) and very inaccurate in finite precision
arithmetic.


\section{Straightforward bounds on the TLS condition number}

\subsection{Sharp lower and upper bounds based on SVD of $[A \,\,b]$ }

In this subsection, we further improve our result in Theorem
\ref{myformationforcond} from the viewpoint of computational cost.
We will show that with the SVD
$$
U^T[A \,\,b] V = {\rm{diag}}(\sigma_1, \ldots, \sigma_{n+1})
$$
we are capable of estimating $\kappa_{g}(A,b)$ accurately based on
the singular values of $[A \,\,b]$ and the last row of $V$ without
calculating $\left\| V^{-T}_{11} S \right\|$, where $V_{11} =
V(1:n,1:n)$,
 $S = {\rm{diag}}([s_1, \ldots, s_n])$, $s_i =
 \frac{\sqrt{\sigma^2_i + \sigma^2_{n+1}}}
 {\sigma^2_i - \sigma^2_{n+1}}$, $i = 1, \ldots,n$,
 as defined in Theorem \ref{myformationforcond}.

From now on we denote $\alpha = \frac{1}{\sqrt{1 + \|x\|^2}}$, which
is always smaller than one for $x\not=0$. Keep
(\ref{V11singularvalues}) in mind and note that
 \begin{equation}\nonumber
  s_1 \leq s_2 \leq \cdots \leq s_n.
 \end{equation}
We then get
 \begin{equation}\nonumber
  s_n  = \sigma_n(V^{-T}_{11}) \|S\| \leq\|V^{-T}_{11} S\| \leq
  \|V^{-T}_{11}\| \|S\| = \alpha^{-1} s_n.
 \end{equation}
Therefore, from Theorem \ref{myformationforcond} we get
 \begin{equation}\label{notmyloveestimate}
\underline{\kappa}:= \alpha^{-1} s_n  \leq \kappa_{g}(A,b) \leq
\bar{\kappa}:=  \alpha^{-2} s_n.
 \end{equation}
So, if $\alpha \approx 1$, that is,  $V_{11}$ is nearly an
orthogonal matrix, the lower and upper bounds in
(\ref{notmyloveestimate}) must be tight.

More generally, for $\alpha$ not small, say, $\frac{1}{2}<\alpha<1$,
we have $\bar{\kappa}< 4 s_n$ and $\underline{\kappa}> s_n$. So
$\underline{\kappa}<\bar{\kappa}<4\underline{\kappa}$. Therefore, in
this case, our lower and upper bounds on the condition number
$\kappa_{g}(A,b)$ are very tight and can estimate the condition
number accurately by no more than four times.

In the following, we only need to discuss the case that
$\alpha\leq\frac{1}{2}$. It will appear that we can establish some
lower bound $\underline{\kappa}$ and upper bound $\bar{\kappa}$ such
that $\underline{\kappa}<\bar{\kappa}<4\underline{\kappa}$ still
holds. As a result, together with the above, for any $0<\alpha<1$,
we can estimate $\kappa_{g}(A,b)$ accurately.

\begin{lemma}\label{CSdecomp}
$V$ can be written as
\begin{equation}\nonumber
V = \left[
            \begin{array}{cc}
              V_{11} & \sqrt{1 - \alpha^2}~ \bar{u}_n \\
              \sqrt{1 - \alpha^2}~ \bar{v}^T_n & - \alpha \\
            \end{array}
          \right],
\end{equation}
where $\bar{u}_n$ and $\bar{v}_n$ are the left and right singular
vectors associated with the smallest singular value of $V_{11}$.
\end{lemma}
{\bf{Proof.}} Based on Theorem \ref{svforV11}, we let
$$
V_{11} = \bar{U} \left[
                        \begin{array}{cc}
                          I_{n-1} &  \\
                           & \alpha \\
                        \end{array}
                      \right] \bar{V}^T
$$
be the SVD of $V_{11}$, where $\bar{U} = [\bar{u}_1, \ldots,
\bar{u}_n] \in \mathbb{R}^{n \times n}$, $\bar{V} =[\bar{v}_1,
\ldots, \bar{v}_n] \in \mathbb{R}^{n \times n}$, and  $\bar{U}^T
\bar{U} = \bar{V}^T \bar{V} = I_n$. It is easily justified from
(\ref{SVDSTLSsol}) that $|V(n+1,n+1)|=\alpha$. Without loss of
generality, we assume $V(n+1, n+1) =-\alpha$. Then, by the theorem
in Section 4 of \cite{PaigeSaunders:1981}, we get
\begin{eqnarray}\nonumber
V &=& \left[
            \begin{array}{cc}
              \bar{U} &  \\
                      & 1 \\
            \end{array}
          \right] \left[
                    \begin{array}{ccc}
                      I_{n-1} & O_{n-1,1} & O_{n-1,1} \\
                     O_{1,n-1}  & \alpha & \sqrt{1-\alpha^2} \\
                      O_{1,n-1} & \sqrt{1-\alpha^2} & -\alpha \\
                    \end{array}
                  \right]
\left[
                           \begin{array}{cc}
                             \bar{V}^T &  \\
                              & 1 \\
                           \end{array}
                         \right]\\ \nonumber 
&=& \left[
      \begin{array}{cc}
        \bar{U}\left[
                 \begin{array}{cc}
                   I_{n-1} &  \\
                    & \alpha \\
                 \end{array}
               \right] \bar{V}^T
         &~~~~ \sqrt{1-\alpha^2} \bar{u}_n  \\
        \sqrt{1-\alpha^2} \bar{v}^T_n  & -\alpha \\
      \end{array}
    \right],
\end{eqnarray}
the desired form of $V$.
 \hfill
$\Box$

Following Lemma \ref{CSdecomp} and letting $[\beta_1, \ldots,
\beta_n, -\alpha]$ be the last row of $V$, we have
\begin{equation}\label{lastrowV}
\bar{v}^T_n = \frac{1}{\sqrt{1 - \alpha^2}} [\beta_1, \ldots,
\beta_n].
\end{equation}
Noticing that ($\alpha^{-1}$, $\bar{u}_n$, $\bar{v}_n$) is the
largest singular triplet of $V^{-T}_{11}$, we denote by
$$
V^{-T}_{11} = [\bar{u}_n, \bar{u}_1 \ldots, \bar{u}_{n-1}] \left[
                                               \begin{array}{cccc}
                                                 \alpha^{-1} &  &  &  \\
                                                  & 1 &  &  \\
                                                  &  & \ddots &  \\
                                                  &  &  & 1 \\
                                               \end{array}
                                             \right] \left[
                                                       \begin{array}{c}
                                                         \bar{v}^T_n \\
                                                         \bar{v}^T_1 \\
                                                         \vdots \\
                                                         \bar{v}^T_{n-1}
                                                       \end{array}
                                                     \right],
$$
which is the SVD of $V^{-T}_{11}$. Then, by (\ref{lastrowV}) we have
\begin{eqnarray}\nonumber
V^{-T}_{11}&=& \left[\alpha^{-1} \bar{u}_n, \bar{u}_1, \ldots,
\bar{u}_{n-1} \right] \left[
\begin{array}{ccccc}
\frac{\beta_1}{\sqrt{1 - \alpha^2}} & \cdots & \frac{\beta_k}
{\sqrt{1 - \alpha^2}} & \cdots & \frac{\beta_n}{\sqrt{1 - \alpha^2}} \\
\bar{v}_1(1) & \ldots & \bar{v}_1(k) & \cdots & \bar{v}_1(n) \\
\vdots &  & \vdots &  & \vdots \\
\bar{v}_{n-1}(1) & \cdots & \bar{v}_{n-1}(k) & \cdots & \bar{v}_{n-1}(n)\\
\end{array}
\right] \\
\label{inV11struct} &=& \left[ \frac{\alpha^{-1} \beta_1}{\sqrt{1 -
\alpha^2}}\bar{u}_n + w_1, \ldots, \frac{\alpha^{-1}
\beta_k}{\sqrt{1 - \alpha^2}}\bar{u}_n + w_k, \ldots,
\frac{\alpha^{-1} \beta_n}{\sqrt{1 - \alpha^2}}\bar{u}_n + w_n
\right],
\end{eqnarray}where $\bar{v}_i(k)$ denotes the $k$th component
of $\bar{v}_i$, $w_k = \sum^{n-1}_{i =1} \bar{v}_i(k)\bar{u}_i$, $k
= 1, \ldots, n$.

\begin{lemma}\label{mylovelemma}
For given matrices $A_1, A_2 \in \mathbb{R}^{n \times n}$, if $A^T_1
A_2 = O$, then
\begin{equation}
\frac{1}{2} (\|A_1\| + \|A_2\|) \leq \|A_1 + A_2\|.
\end{equation}
\end{lemma}
{\bf{Proof.}} For an arbitrary vector $x \in \mathbb{R}^{n}$, from
$(A_1 x)^T (A_2 x) = 0$ it follows that
\begin{equation}\nonumber
\|A_1 x\|, \|A_2 x\| \leq \|A_1 x + A_2 x\|
\end{equation}and that
\begin{eqnarray}\nonumber
\|A_1\| &=& {\rm{max}}_{\|x\| =1} \|A_1 x\| \leq {\rm{max}}_{\|x\|
=1} \|A_1 x + A_2 x\| = \|A_1 + A_2\|,
\\\nonumber
\|A_2\| &=& {\rm{max}}_{\|x\| =1} \|A_2 x\| \leq {\rm{max}}_{\|x\|
=1} \|A_1 x + A_2 x\| = \|A_1 + A_2\|.
\end{eqnarray}
So, we get the desired inequality. \hfill $\Box$

To prove the main results of this section, we need the following two
propositions.
\begin{proposition}\label{bigsurpriseprop}
Let $[\beta_1, \ldots, \beta_n, -\alpha]$ be the last row of $V$,
$V_{11}=V(1:n, 1:n)$ and $\bar{S} = {\rm{diag}}([\bar{s}_1, \ldots,
\bar{s}_n])$, where $\bar{s}_1, \ldots, \bar{s}_n$ are arbitrary
positive numbers and satisfy $0< \bar{s}_1 \leq \bar{s}_2 \leq
\cdots \leq \bar{s}_n$. Then
\begin{eqnarray}\nonumber 
&& \underline{c} := \frac{1}{2} \left( \frac{\alpha^{-1}
\sqrt{\beta^2_1 \bar{s}^2_1 + \ldots + \beta^2_n \bar{s}^2_n }
}{\sqrt{1 - \alpha^{2}}} + \frac{\sqrt{1 - \alpha^2 - \beta^2_n  }}{
\sqrt{1 - \alpha^2} }  \bar{s}_n \right) \\ \nonumber &\leq &
\left\|V^{-T}_{11} \bar{S} \right\|\leq \bar{c} := \frac{\alpha^{-1}
\sqrt{\beta^2_1 \bar{s}^2_1 + \ldots + \beta^2_n \bar{s}^2_n }
}{\sqrt{1 - \alpha^{2}}} + \bar{s}_n.
\end{eqnarray}
\end{proposition}
{\bf{Proof.}} Following (\ref{inV11struct}), we get
\begin{eqnarray}\nonumber 
&&V^{-T}_{11} \bar{S} \\ \nonumber &=& \left[ \frac{\alpha^{-1}
\beta_1 \bar{s}_1}{\sqrt{1 - \alpha^2}}\bar{u}_n + \bar{s}_1 w_1,
\ldots, \frac{\alpha^{-1} \beta_k \bar{s}_k}{\sqrt{1 -
\alpha^2}}\bar{u}_n + \bar{s}_k w_k, \ldots, \frac{\alpha^{-1}
\beta_n \bar{s}_n}{\sqrt{1 - \alpha^2}}\bar{u}_n + \bar{s}_n w_n
\right].
\end{eqnarray}
Define
\begin{equation}\nonumber
A_1 = \left[\frac{\alpha^{-1} \beta_1 \bar{s}_1}{\sqrt{1 -
\alpha^2}}\bar{u}_n, \ldots, \frac{\alpha^{-1} \beta_n
\bar{s}_n}{\sqrt{1 - \alpha^2}}\bar{u}_n \right], \,\, A_2 =
\left[\bar{s}_1 w_1, \ldots, \bar{s}_n w_n   \right].
\end{equation}Then $V^{-T}_{11} \bar{S} = A_1 + A_2$. Noticing that
\begin{equation}\nonumber
\bar{u}^T_n w_k = 0,\,\, k = 1, \ldots, n,
\end{equation}
we get $A^T_1 A_2 = O$. Thus, we have
\begin{equation}\label{surprise1}
\frac{1}{2} (\|A_1\| + \|A_2\|) \leq \left\|V^{-T}_{11} \bar{S}
\right\| \leq \|A_1\| + \|A_2\|,
\end{equation}
in which the left-hand side inequality follows from Lemma
\ref{mylovelemma}. Furthermore, noticing that
\begin{equation}\nonumber
A_1 = \frac{\alpha^{-1}}{\sqrt{1 - \alpha^2}}\bar{u}_n \left[\beta_1
\bar{s}_1, \ldots,  \beta_n \bar{s}_n\right]
\end{equation}
and $\|\bar u_n\|=1$, we have
\begin{equation}\label{surprise2}
\|A_1\| = \frac{\alpha^{-1}}{\sqrt{1 - \alpha^2}} \left\|
\left[\beta_1 \bar{s}_1, \ldots, \beta_n \bar{s}_n \right] \right\|
= \frac{\alpha^{-1}}{\sqrt{1 - \alpha^2}} \sqrt{\beta^2_1
\bar{s}^2_1+ \ldots+ \beta^2_n \bar{s}^2_n}.
\end{equation} In the meantime, note that
\begin{equation}\nonumber
\|w_n\| = \sqrt{\sum^{n-1}_{i=1} \bar{v}^2_i(n)}=\sqrt{1-
\frac{\beta^2_n}{1-\alpha^2}} =
\frac{\sqrt{1-\alpha^2-\beta^2_n}}{\sqrt{1-\alpha^2}},
\end{equation}
\begin{equation}\nonumber
\|[w_1, \ldots, w_n]\| =\left\| [\bar{u}_1, \ldots, \bar{u}_{n-1}]
\left[
                                                          \begin{array}{c}
                                                            \bar{v}^T_1 \\
                                                            \vdots \\
                                                            \bar{v}^T_{n-1} \\
                                                          \end{array}
                                                        \right] \right\| =1,
\end{equation}and
\begin{equation}\nonumber
\|\bar{S}\| = \bar{s}_n.
\end{equation}
From
\begin{equation}\nonumber
\|\bar{s}_n w_n\| \leq \|A_2\| \leq \left\|[w_1, \ldots, w_n]
\right\| \|\bar{S}\|
\end{equation}
we get
\begin{equation}\label{surprise3}
\frac{\sqrt{1-\alpha^2 - \beta^2_n}}{\sqrt{1-\alpha^2}}\bar{s}_n
\leq \|A_2\| \leq \bar{s}_n.
\end{equation}
Combining (\ref{surprise1}), (\ref{surprise2}) and
(\ref{surprise3}), we establish the desired inequality. \hfill
$\Box$

\begin{proposition}\label{oneremark}
Suppose that $\alpha \leq\frac{1}{2}$. Then for $\underline{c}$ and
$\bar{c}$ in Proposition \ref{bigsurpriseprop}, we have
\begin{equation}\label{remarkineqn}
\underline{c} < \bar{c}< 4 \underline{c}.
\end{equation}
\end{proposition}
{\bf{Proof.}} If
$\frac{|\beta_n|}{\sqrt{1-\alpha^2}}<\frac{\sqrt{3}}{2}$, then it is
easy to verify that
\begin{equation}\nonumber
\frac{\sqrt{1-\alpha^2 - \beta^2_n}}{\sqrt{1-\alpha^2}}
>\frac{1}{2}
\end{equation}
and
\begin{equation}\nonumber
        \underline{c} >\frac{1}{4} \bar{c}.
\end{equation}
Thus, (\ref{remarkineqn}) holds. If
$\frac{|\beta_n|}{\sqrt{1-\alpha^2}}\geq\frac{\sqrt{3}}{2}$, then
$$
\alpha^{-1} \frac{|\beta_n|}{\sqrt{1-\alpha^2}}
>\frac{\sqrt{3}}{2}\alpha^{-1} >1,
$$
so $\alpha^{-1} \frac{|\beta_n|}{\sqrt{1-\alpha^2}} \bar{s}_n >
\bar{s}_n$, from which and the definitions of $\bar{c}$ and
$\underline{c}$ it follows that
\begin{eqnarray}\nonumber
\bar{c} &< & \frac{\alpha^{-1} \sqrt{\beta^2_1 \bar{s}^2_1 + \ldots
+ \beta^2_n \bar{s}^2_n }  }{\sqrt{1 - \alpha^{2}}} + \alpha^{-1}
\frac{|\beta_n|}{\sqrt{1-\alpha^2}} \bar{s}_n
\\\nonumber
&\leq & \frac{2 \alpha^{-1} \sqrt{\beta^2_1 \bar{s}^2_1 + \ldots +
\beta^2_n \bar{s}^2_n }  }{\sqrt{1 - \alpha^{2}}}
\\\nonumber
&\leq & \frac{2 \alpha^{-1} \sqrt{\beta^2_1 \bar{s}^2_1 + \ldots +
\beta^2_n \bar{s}^2_n }  }{\sqrt{1 - \alpha^{2}}} +  \frac{2 \sqrt{1
- \alpha^2 - \beta^2_n  }}{ \sqrt{1 - \alpha^2} }  \bar{s}_n = 4
\underline{c}.
\end{eqnarray}
Thus, (\ref{remarkineqn}) still holds. \hfill $\Box$

Now we are in a position to derive sharp bounds on
$\kappa_{g}(A,b)$.
\begin{theorem}\label{bigsurprisethm}
Let $[\beta_1, \ldots, \beta_n, -\alpha]$ be the last row of $V$ and
$S = {\rm{diag}} ([s_1, \ldots, s_n])$, $s_i =
\frac{\sqrt{\sigma^2_i + \sigma^2_{n+1}}}{\sigma^2_i -
\sigma^2_{n+1}}$, $i = 1, \ldots, n$. Then
\begin{eqnarray}\nonumber 
&& \underline{\kappa} := \frac{1}{2} \left( \frac{\alpha^{-2}
\sqrt{\beta^2_1 s^2_1 + \ldots + \beta^2_n s^2_n }  }{\sqrt{1 -
\alpha^{2}}} + \frac{\sqrt{1 - \alpha^2 - \beta^2_n  }}{ \sqrt{1 -
\alpha^2} } \alpha^{-1} s_n \right)
\\\nonumber
&\leq & \kappa_{g}(A,b)\leq \bar{\kappa} := \frac{\alpha^{-2}
\sqrt{\beta^2_1 s^2_1 + \ldots + \beta^2_n s^2_n }  } {\sqrt{1 -
\alpha^{2}}} + \alpha^{-1}s_n.
\end{eqnarray}Moreover, if $\alpha \leq \frac{1}{2}$, then
\begin{equation}\nonumber
\underline{\kappa} < \bar{\kappa} < 4 \underline{\kappa}.
\end{equation}
\end{theorem}
{\bf{Proof.}} Noticing that $0< s_1 \leq s_2 \leq \cdots \leq s_n$
and using Proposition \ref{bigsurpriseprop}, we have
\begin{eqnarray}\nonumber
&& \frac{1}{2} \left( \frac{\alpha^{-1} \sqrt{\beta^2_1 s^2_1 +
\ldots + \beta^2_n s^2_n }  }{\sqrt{1 - \alpha^{2}}} + \frac{\sqrt{1
- \alpha^2 - \beta^2_n  }}{ \sqrt{1 - \alpha^2} } s_n \right) \\
\nonumber 
&\leq & \left\|V^{-T}_{11} S \right\| \leq \frac{\alpha^{-1}
\sqrt{\beta^2_1 s^2_1 + \ldots + \beta^2_n s^2_n } }{\sqrt{1 -
\alpha^{2}}} + s_n.
\end{eqnarray} By Theorem \ref{myformationforcond}, we get the first part of the
theorem. Furthermore, we have the second part of the theorem  by
Proposition \ref{oneremark}. \hfill $\Box$

{\bf A small example (Continued)}. From Theorem
\ref{bigsurprisethm}, we have
\begin{equation}\nonumber
 5.65\times 10^8\leq \kappa^r_{g}(A,b) \leq 1.13\times 10^9.
\end{equation}
The lower and upper bounds estimate
${\kappa}^r_g(\ref{mynewclosed})=1.13\times 10^9$ accurately, as
described in the second part of Theorem \ref{bigsurprisethm}.

\subsection{Lower and upper bounds based on a few of singular values of
$A$ and $[A \,\,b]$}

In \cite{Malyshev:2003}, bounds on the condition number of the
Tikhonov regularization solution are established based on a few
singular values of $A$, where $A$ is the coefficient matrix of the
least squares problem under consideration. This is particularly
useful for large scale TLS problems since for them any formulas and
bounds for the condition number involving all the singular values of
$A$ and/or $[A \ b]$ are too costly to be computed. Such a bound can
be obtained by computing only a few singular values of $A$ and/or
$[A\ b]$.

In the following theorem, we establish similar results for the
condition number of the TLS problem.

\begin{theorem}\label{firstestimateforcond}
We have
\begin{equation}\label{luestimate}
\underline{\kappa}_1 \leq \kappa_{g}(A,b) \leq  \bar{\kappa}_1,
\end{equation}
where
\begin{eqnarray}\label{lubound}
\underline{\kappa}_1 = \frac{\sqrt{1 + \|x\|^2} \sqrt{
\hat{\sigma}^2_{n-1} + \sigma^2_{n+1} }}{ \hat{\sigma}^2_{n-1} -
\sigma^2_{n+1}}, \,\, \bar{\kappa}_1 = \frac{\sqrt{1 + \|x\|^2}
\sqrt{ \hat{\sigma}^2_{n} + \sigma^2_{n+1} }} {  \hat{\sigma}^2_{n}
- \sigma^2_{n+1}}.
\end{eqnarray}
\end{theorem}

{\bf{Proof.}} Denoting
$$M = (A^T A - \sigma^2_{n+1} I_n)^{-1} \left( (
\|x\|^2 + 1) A^T A +  \|r\|^2 I_n \right) (A^T A - \sigma^2_{n+1}
I_n)^{-1},
$$
from (\ref{firstchange}) we have
\begin{equation}\label{fromotherproof}
K  K^T = M - 2 \sigma^2_{n+1} (A^T A - \sigma^2_{n+1} I_n)^{-1} x
x^T (A^T A - \sigma^2_{n+1} I_n)^{-1}.
\end{equation}
Here and hereafter, $\lambda_i(M)$ denotes the $i$th largest
eigenvalue of $M$, where $M$ is an arbitrary symmetric matrix. By
the Courant-Fischer theorem \cite[p.182]{Horn:1999}, from
(\ref{fromotherproof}) we get
\begin{equation}\label{inequality1}
 \lambda_2(M) \leq  \lambda_1(K K^T).
\end{equation}
Furthermore, since $2\sigma^2_{n+1}(A^T A - \sigma^2_{n+1} I_n)^{-1}
x x^T (A^T A - \sigma^2_{n+1} I_n)^{-1}$ is nonnegative definite,
the following inequality holds
\begin{equation}\label{inequlity2}
 \lambda_1(K  K^T) \leq \lambda_1(M).
\end{equation}
Collecting (\ref{inequality1}) and (\ref{inequlity2}) and based on
(\ref{knoneckercond}), we have
\begin{equation}\nonumber
\sqrt{\lambda_2(M)}  \leq \kappa_{g}(A,b) \leq \sqrt{\lambda_1(M)}.
\end{equation}

It is easy to verify that the set
\begin{equation}\nonumber
\left\{ \frac{(\|x\|^2 + 1) \hat{\sigma}^2_j +
 \|r\|^2}{(\hat{\sigma}^2_j - \sigma^2_{n+1})^2}
\right\}^{n}_{j=1}
\end{equation} consists of all the eigenvalues of $M$. We define the function
\begin{equation}\nonumber
f(\sigma) = \frac{( \|x\|^2 + 1) {\sigma}^2 +
 \|r\|^2}{({\sigma}^2 - \sigma^2_{n+1})^2}, \,\,\sigma >
\sigma_{n+1},
\end{equation}and differentiate it to get
\begin{equation}\nonumber
f'(\sigma) = \frac{-2 \sigma^3 (\|x\|^2 + 1) - 2 \sigma ( \|x\|^2 +
1) \sigma^2_{n+1} - 4 \sigma \|r\|^2}{(\sigma^2 -
\sigma^2_{n+1})^3}.
\end{equation}
It is seen that $f'(\sigma) < 0$ and $f(\sigma)$ is decreasing in
the interval $(\sigma_{n+1}, \infty) $. Thus, we get that
\begin{equation}\nonumber
 \lambda_1(M) = \frac{(\|x\|^2 + 1)
\hat{\sigma}^2_n + \|r\|^2}{(\hat{\sigma}^2_n -
\sigma^2_{n+1})^2},\,\, \lambda_2(M) = \frac{(\|x\|^2 + 1)
\hat{\sigma}^2_{n-1} + \|r\|^2}{(\hat{\sigma}^2_{n-1} -
\sigma^2_{n+1})^2}
\end{equation}
and
\begin{equation}\nonumber
\frac{\sqrt{(\|x\|^2 + 1) \hat{\sigma}^2_{n-1} +
\|r\|^2}}{\hat{\sigma}^2_{n-1} - \sigma^2_{n+1}} \leq
\kappa_{g}(A,b) \leq
    \frac{\sqrt{( \|x\|^2 + 1) \hat{\sigma}^2_n +
\|r\|^2}}{\hat{\sigma}^2_n - \sigma^2_{n+1}}.
\end{equation}

Noticing that $\frac{\|r\|^2}{1 + \|x\|^2} = \sigma^2_{n+1}$, we
complete the proof. \hfill $\Box$

{\bf Remark}. In Corollary 1 of \cite{BaboulinGratton:2010}, the
authors prove that
\begin{equation}\nonumber
\kappa_{g}(A,b) \leq \frac{\sqrt{1 + \|x\|^2} \sqrt{ {\sigma}^2_{1}
+\sigma^2_{n+1} }} {  \hat{\sigma}^2_{n} - \sigma^2_{n+1}}.
\end{equation}
Since $\hat{\sigma}_n\leq\hat{\sigma}_1,\
\hat{\sigma}_1\leq\sigma_1$, we get
\begin{equation}\nonumber
\bar{\kappa}_1 \leq \frac{\sqrt{1 + \|x\|^2} \sqrt{
\hat{\sigma}^2_{1} + \sigma^2_{n+1} }}{  \hat{\sigma}^2_{n} -
\sigma^2_{n+1}} \leq \frac{\sqrt{1 + \|x\|^2} \sqrt{ {\sigma}^2_{1}
+\sigma^2_{n+1} }} {\hat{\sigma}^2_{n} - \sigma^2_{n+1}}.
\end{equation}
Therefore, our $\bar{\kappa}_1$ in (\ref{lubound}) is sharper than
the above upper bound.

It is seen that the lower and upper bounds on $\kappa_{g}(A,b)$ in
Theorem \ref{firstestimateforcond} are marginally different provided
that  $\hat{\sigma}_n$ and $\hat{\sigma}_{n-1}$ are close. This
means that in this case both bounds are very tight. For the case
that $\hat{\sigma}_n$ and $\hat{\sigma}_{n-1}$ are not close, we
next give a new lower bound that can be better than that in Theorem
\ref{firstestimateforcond}.

\begin{theorem}\label{myfirstestimate}
It holds that
\begin{equation}\nonumber
\underline{\kappa}_2 \leq \kappa_{g}(A,b)\leq \bar{\kappa}_1,
\end{equation}where $\bar{\kappa}_1$ is defined as in
Theorem \ref{firstestimateforcond} and
\begin{equation}\nonumber
\underline{\kappa}_2 =\frac{\sqrt{1 + \|x\|^2}}{  \sqrt{
{\hat{\sigma}}^2_n
 - \sigma^2_{n+1}}}.
\end{equation}Moreover, when $\hat{\sigma}_{n-1} \geq \sigma_{n+1} +
\sqrt{{\hat{\sigma}}^2_n - \sigma^2_{n+1}}$, we have
\begin{equation}\nonumber
\underline{\kappa}_1 \leq \underline{\kappa}_2. \end{equation}
\end{theorem}
{\bf{Proof.}} Denote $P = A^T A - \sigma^2_{n+1}I_n$. From
(\ref{forcond}), we have
\begin{equation}\nonumber
 \frac{1}{\|x\|^2 + 1} ~K  K^T =P^{-1} + 2\sigma^2_{n+1} P^{-1}
 \left( I_n - \frac{ x x^T}{1 +
\|x\|^2} \right) P^{-1}.
\end{equation}
Noticing the second term in the right-hand side of the above
relation is positive definite, we have
\begin{equation}\nonumber
( \|x\|^2 +1 )\lambda_1(P^{-1}) \leq \lambda_1(  ~K K^T ),
\end{equation}that is,
\begin{equation}\nonumber
\frac{ \|x\|^2 + 1}{{\hat{\sigma}^2_n} - \sigma^2_{n+1}}\leq
\kappa^2_{g}(A,b).
\end{equation} Thus, the first part of the theorem is obtained.

The second part of the theorem is proved by noting
\begin{equation}\nonumber
\frac{ \sqrt{ \hat{\sigma}^2_{n-1} + \sigma^2_{n+1} }}{
\hat{\sigma}^2_{n-1} - \sigma^2_{n+1}}< \frac{1}{ \hat{\sigma}_{n-1}
- \sigma_{n+1}}\leq \frac{1} {\sqrt{{\hat{\sigma}}^2_n -
\sigma^2_{n+1}}}
\end{equation}under the assumption that $\hat{\sigma}_{n-1} -\sigma_{n+1}\geq
\sqrt{{\hat{\sigma}}^2_n - \sigma^2_{n+1}}$. \hfill $\Box$

{\bf{Remark 1.}} At first glance, the assumption in the second part
of the theorem seems not so direct but we can justify that it indeed
implies that $\hat{\sigma}_n$ and $\hat{\sigma}_{n-1}$ are not
close. Actually, we can verify that the second part of Theorem
\ref{myfirstestimate} holds under a slightly stronger but much
simpler condition that
\begin{equation}\nonumber
 \hat{\sigma}_{n-1} \geq 2 \hat{\sigma}_n.
 \end{equation}

{\bf{Remark 2.}} From
\begin{equation}\nonumber
\frac{\bar{\kappa}_1}{\underline{\kappa}_2} =
\frac{\sqrt{\hat{\sigma}^2_n +  \sigma^2_{n+1}}}
{\sqrt{\hat{\sigma}^2_n - \sigma^2_{n+1}}} = \sqrt{\frac{1 +
\frac{\sigma^2_{n+1}}{\hat{\sigma}^2_n}} {1 -  \frac{\sigma^2_{n+1}}
{\hat{\sigma}^2_n}  }},
\end{equation}it is seen that $\frac{\bar{\kappa}_1}{\underline{\kappa}_2}>
1$  provided $\sigma_{n+1} > 0$. Only for $\sigma_{n+1} = 0$,
$\bar{\kappa}_1 = \underline{\kappa}_2$ holds. At this time, $b \in
\mathcal{R}(A)$ and $r = 0$.

We observe that the bounds on $\kappa_g(A,b)$ in Theorem
\ref{myfirstestimate} are tight when
$\frac{\sigma_{n+1}}{\hat{\sigma}_n}$ is small, compared with one.
On the other hand, once $\frac{\sigma_{n+1}}{\hat{\sigma}_n} $ is
not small, these bounds may not be tight. For this case, we will
present new bounds that may better estimate $\kappa_g(A,b)$.

The proof of the following theorem depends strongly on Propositions
\ref{bigsurpriseprop} and \ref{oneremark}.

\begin{theorem}\label{improveTLSthm}
Assume that $\alpha \leq\frac{1}{2}$. Denote $\rho =
\frac{\sigma_{n+1}}{\sigma_n}$. Then
\begin{eqnarray}\label{TLSimprove}
\underline{\kappa}_2: =\frac{\sqrt{1 +
\|x\|^2}}{\sqrt{\hat{\sigma}^2_n - \sigma^2_{n+1}}} \leq
 \kappa_g(A,b) &<& \bar{\kappa}_2:=\sqrt{\frac{1+31 \rho^2}{1-\rho^2}}
 \frac{\sqrt{1+ \|x\|^2}
}{\sqrt{\hat{\sigma}^2_n - \sigma^2_{n+1}}}.
\end{eqnarray}
\end{theorem}

{\bf{Proof.}} Based on Theorem \ref{myfirstestimate}, it suffices to
prove the right-hand side of (\ref{TLSimprove}). From
(\ref{forcond}) and (\ref{Wewilluseit}), we get
\begin{eqnarray}\nonumber
 \frac{1}{\|x\|^2 + 1 }~K K^T &=& P^{-1} + 2 \sigma^2_{n+1} P^{-1}
 \left( I_n  -  \frac{ x x^T}{1 + \|x\|^2} \right) P^{-1}, \\ \label{fornewbound}
 &=& V^{-T}_{11} \Lambda^{-1} V^{-1}_{11} + 2 \sigma^2_{n+1} V^{-T}_{11}
 \Lambda^{-2} V^{-1}_{11}:=P^{-1} + E,
\end{eqnarray}
where $P=A^T A - \sigma^2_{n+1} I_n$, $ \Lambda = {\rm{diag}}
([\sigma^2_1 - \sigma^2_{n+1}, \ldots, \sigma^2_n -
\sigma^2_{n+1}])$. Denote
\begin{eqnarray}\nonumber
D &=& {\rm{diag}}([d_1, \ldots, d_n]),\, d_i = \frac{\sigma_{n+1}}
{\sigma^2_i - \sigma^2_{n+1}}, i = 1, \ldots, n,
\\\nonumber
T &=& {\rm{diag}}([t_1, \ldots, t_n]),\, t_i =
\frac{1}{\sqrt{\sigma^2_i - \sigma^2_{n+1}}}, i = 1, \ldots, n.
\end{eqnarray}
Then $P^{-1} = \left( V^{-T}_{11} T\right) \left( T V^{-1}_{11}
\right)$ and $E = 2 \left( V^{-T}_{11} D\right)  \left(D V^{-1}_{11}
\right)$.

Note that $0 < d_1 \leq d_2 \leq \cdots \leq d_n$ and $0 < t_1 \leq
t_2 \leq \cdots \leq t_n$. Applying Proposition
\ref{bigsurpriseprop}, we get
\begin{eqnarray}\nonumber 
&& \frac{1}{2} \left( \frac{\alpha^{-1} \sqrt{\beta^2_1 d^2_1 +
\ldots + \beta^2_n d^2_n }  }{\sqrt{1 - \alpha^{2}}} + \frac{\sqrt{1
- \alpha^2 - \beta^2_n  }}{ \sqrt{1 - \alpha^2} } d_n \right) \\
\label{lambdabigsurprise3} &\leq & \left\|V^{-T}_{11} D  \right\|
\leq \frac{\alpha^{-1} \sqrt{\beta^2_1 d^2_1 + \ldots + \beta^2_n
d^2_n }  }{\sqrt{1 - \alpha^{2}}} + d_n
\end{eqnarray}
and
\begin{eqnarray}\nonumber
&& \frac{1}{2} \left( \frac{\alpha^{-1} \sqrt{\beta^2_1 t^2_1 +
\ldots + \beta^2_n t^2_n }  }{\sqrt{1 - \alpha^{2}}} + \frac{\sqrt{1
- \alpha^2 - \beta^2_n  }}{ \sqrt{1 - \alpha^2} } t_n \right) \\
\nonumber
 &\leq & \left\|V^{-T}_{11} T  \right\|
\leq \frac{\alpha^{-1} \sqrt{\beta^2_1 t^2_1 + \ldots + \beta^2_n
t^2_n }  }{\sqrt{1 - \alpha^{2}}} + t_n,
\end{eqnarray}respectively, where $[\beta_1, \ldots,
\beta_n, -\alpha]$ denotes the last row of $V$ as before. Define
$k_n = \frac{d_n}{t_n}=\frac{\sigma_{n+1}}{\sqrt{\sigma^2_n -
\sigma^2_{n+1}   }}$. Then
\begin{equation}\nonumber
\frac{d_1}{t_1}= \frac{\sigma_{n+1}}{\sqrt{\sigma^2_1 -
\sigma^2_{n+1}   }} \leq k_n\,, \ldots, \frac{d_{n-1}}{t_{n-1}}
=\frac{\sigma_{n+1}}{\sqrt{\sigma^2_{n-1} - \sigma^2_{n+1}   }} \leq
k_n.
\end{equation}
Thus, by (\ref{lambdabigsurprise3}) we have
\begin{eqnarray}\label{firstimproineqn}
\frac{1}{\sqrt{2}} \|E\|^{\frac{1}{2}} = \left\|V^{-T}_{11} D
\right\|   \leq  k_n \left( \frac{\alpha^{-1}
 \sqrt{\beta^2_1 t^2_1 + \ldots + \beta^2_n t^2_n }  }{\sqrt{1 - \alpha^{2}}} + t_n
 \right).
\end{eqnarray}
Note that for the lower and upper bounds on $\left\|V^{-T}_{11} T
\right\|$ above, by Proposition \ref{oneremark} it holds that
\begin{eqnarray}\nonumber
\frac{\alpha^{-1} \sqrt{\beta^2_1 t^2_1 + \ldots + \beta^2_n t^2_n }
}{\sqrt{1 - \alpha^{2}}} + t_n &<& 2 \left( \frac{\alpha^{-1}
\sqrt{\beta^2_1 t^2_1 + \ldots + \beta^2_n t^2_n } }{\sqrt{1 -
\alpha^{2}}} + \frac{\sqrt{1 - \alpha^2 - \beta^2_n  }}{ \sqrt{1 -
\alpha^2} } t_n \right) \\ \label{lambdabigsurprise4} &<& 4
\left\|V^{-T}_{11} T \right\|.
\end{eqnarray}
Based on (\ref{firstimproineqn}) and (\ref{lambdabigsurprise4}), we
derive that
\begin{equation}\nonumber
\frac{1}{\sqrt{2}} \|E\|^{\frac{1}{2}} < 4 k_n \left\|V^{-T}_{11} T
\right\| = 4 k_n \|P^{-1}\|^{\frac{1}{2}}
\end{equation}
and that
\begin{equation}\label{improveE}
\|E\| < 32 k^2_n \|P^{-1}\|.
\end{equation}
Combining  (\ref{improveE}) and (\ref{fornewbound}), we establish
that
\begin{eqnarray}\nonumber
\kappa_g(A,b) = \|K\| = \|K K^T\|^{\frac{1}{2}} & < & \sqrt{1 +
32k^2_n} \sqrt{1+ \|x\|^2} \|P^{-1}\|^{\frac{1}{2}}\\\nonumber 
&=& \sqrt{\frac{1+31 \rho^2}{1-\rho^2}    } \frac{\sqrt{1+
\|x\|^2}}{ \sqrt{\hat{\sigma}^2_n - \sigma^2_{n+1} } }.
\end{eqnarray}So, the proof of the theorem is completed. \hfill $\Box$

{\bf{Remark.}} It is clear that the bounds in Theorem
\ref{improveTLSthm} are tight when
$\rho=\frac{\sigma_{n+1}}{\sigma_{n}}$ is small, compared with one.
The result in this theorem is of particular importance in the case
that $\frac{\sigma_{n+1}}{\hat{\sigma}_n} $ is close to one. Recall
that the lower and upper bounds in Theorem \ref{myfirstestimate}
differ considerably when $\frac{\sigma_{n+1}}{\hat{\sigma}_n} $ is
close to one. Theorem \ref{improveTLSthm} tells us that, if only
$\frac{\sigma_{n+1}}{\sigma_{n}}$ is not so close to one,
$\kappa_g(A,b)$ should be close to the lower bound.

The improvement of $\bar{\kappa}_2$ to $\bar{\kappa}_1$ can be
illustrated as follows. For $\frac{\sigma_{n+1}}{\sigma_{n}}$ small,
i.e., ${\sigma_{n+1}}$ and ${\sigma_{n}}$ not close, as an upper
bound of $\kappa^r_g(A,b)$,
\begin{eqnarray}\nonumber
\bar{\kappa}^r_2 := \frac{\bar{\kappa}_2}{\|x\|} \|[A\,\,b]\|_F &=&
\sqrt{\frac{1+31 \rho^2}{1-\rho^2}}\frac{\sqrt{1+ \|x\|^2}}{\|x\|}
\frac{\|[A\,\,b]\|_F}{\sqrt{\hat{\sigma}^2_n - \sigma^2_{n+1} }}
\\\nonumber
&\approx & \sqrt{\frac{1+31\rho^2}{1-\rho^2}}
\frac{\|[A\,\,b]\|_F}{\sqrt{\hat{\sigma}^2_n - \sigma^2_{n+1} }}
\end{eqnarray} is a moderate multiple of $\frac{1}
{\sqrt{\hat{\sigma}^2_n - \sigma^2_{n+1} }}$. In contrast,
\begin{eqnarray}\nonumber
\bar{\kappa}^r_1 := \frac{\bar{\kappa}_1}{\|x\|} \|[A\,\,b]\|_F &=&
\frac{\sqrt{1+ \|x\|^2}}{\|x\|} \frac{ \sqrt{\hat{\sigma}^2_n +
\sigma^2_{n+1} }    }{ \hat{\sigma}^2_n - \sigma^2_{n+1} }
\|[A\,\,b]\|_F \\ \nonumber &\approx & \frac{ \sqrt{\hat{\sigma}^2_n
+ \sigma^2_{n+1} }    }{ \hat{\sigma}^2_n - \sigma^2_{n+1} }
\|[A\,\,b]\|_F
\end{eqnarray}
is a moderate multiple of $\frac{1} {\hat{\sigma}^2_n -
\sigma^2_{n+1} }$. The improvement of $\bar{\kappa}^r_2$ over
$\bar{\kappa}^r_1$ becomes significant as ${\sigma_{n+1}}$ and
${\hat{\sigma}_n} $ are close. Similarly, $\bar{\kappa}^r_2$ also
improves the approximate condition number used in
\cite{BjorckHeggernesMatstoms:2000}:
\begin{equation}\nonumber
\bar{\kappa}^r_{\cite{BjorckHeggernesMatstoms:2000}} :=
\frac{\hat{\sigma}_1 }{\hat{\sigma}_n - \sigma_{n+1}    }  = \frac{
\hat{\sigma}_1 (\hat{\sigma}_n + \sigma_{n+1})    }{\hat{\sigma}^2_n
- \sigma^2_{n+1}}.
\end{equation}

We will further illustrate the improvement by numerical experiments
to be presented in Section 5.

\section{Numerical experiments}

We present numerical experiments to illustrate how tight the bounds
in Theorems \ref{myfirstestimate} and \ref{improveTLSthm} are, and
to compare the bounds with the related result in
\cite{BjorckHeggernesMatstoms:2000}. For a given TLS problem, the
TLS solution is computed by (\ref{SVDSTLSsol}). All experiments were
run using Matlab 7.8.0 with the machine precision $\epsilon_{\rm
mach}=2.22\times 10^{-16}$ under the Microsoft Windows XP operating
system.

{\bf{Example 1.}} In this example, the TLS problem comes from
\cite{KammNagy:1998}. Specifically, an $m \times (m - 2 \omega) $
convolution matrix $\bar{T}$ is constructed to have the first column
\begin{equation}\nonumber
t_{i,1} = \left\{
  \begin{array}{ll}
    \frac{1}{\sqrt{2 \pi \alpha^2} {\rm{exp}}
    \left[ \frac{-(\omega -i+1)^2}{2 \alpha^2} \right]} &
    \hbox{\,\,\,$i = 1, 2, \ldots, 2 \omega  +1$,} \\
    0 & \hbox{\,\,\,otherwise,}
  \end{array}
\right.
\end{equation}
and the first row
\begin{equation}\nonumber
t_{1,j} = \left\{
            \begin{array}{ll}
              t_{1,1} & \hbox{\,\,\,if $j =1$,} \\
              0 & \hbox{\,\,\,otherwise,}
            \end{array}
          \right.
\end{equation} where $\alpha =
1.25$ and $\omega = 8$.  A Toeplitz matrix $A$ and a right-hand side
vector $b$ are constructed as $A = \bar{T} + E$ and $b = \bar{g} +
e$, where $\bar{g} =[1, \ldots, 1]^T$, $E$ is a random Toeplitz
matrix with the same structure as $\bar{T}$ and $e$ is a random
vector. The entries in $E$ and $e$ are generated randomly from a
normal distribution with mean zero and variance one, and scaled so
that
\begin{equation}\nonumber
\|e\| = \gamma \|\bar{g}\|, \,\,\,\|E\| = \gamma \|\bar{T}\|,
\,\,\gamma = 0.001.
\end{equation}

\begin{table}[h]\label{tab4}
\begin{center}
\caption{}
{ \doublerulesep18.0pt \tabcolsep0.01in
\begin{tabular}{|c |c| c | c| c | c |c |c|}
\hline $m$   & $\sigma_{n+1}/\sigma_n$ &
$\sigma_{n+1}/\hat{\sigma}_n$& $\kappa^r_g(A,b)$ & 
 $\underline{\kappa}^r_{2}$& $\bar{\kappa}^r_{2} $  & $\bar{\kappa}^r_{1} $  &
 $\bar{\kappa}^r_{\cite{BjorckHeggernesMatstoms:2000}} $
 \\
\hline
 $~100~$   & $0.981$ &
 $~ 1-7.85\times 10^{-9}$ & $7.70\times 10^7$   & $7.04\times 10^7$&
 $2.01\times 10^9$ &$7.94\times 10^{11}$
 &$1.03\times 10^{11}$
 \\
 \hline
  $~300~$   & $0.995$ &
 $~ 1-2.05\times 10^{-8}$ & $1.40\times 10^8$
 & $1.26\times 10^8$& $6.90\times 10^9$ &$8.83\times 10^{11}$
 &$6.54\times 10^{10}$
 \\
 \hline
  $~500~$   & $0.998$ &
 $~ 1-5.66\times 10^{-8}$ & $9.01\times 10^7$   & $7.89\times 10^7$
 & $6.56\times 10^9$ &$3.32\times 10^{11}$
 &$1.90\times 10^{10}$
 \\
 \hline
\end{tabular}}
\end{center}
\end{table}

In the table,
$$
\underline{\kappa}^r_2 = \frac{\underline{\kappa}_2 }{\|x\|}\|[A
\,\,b]\|_F, \,\, \bar{\kappa}^r_2 = \frac{\bar{\kappa}_2 }{\|x\|}
\|[A \,\,b]\|_F, \,\, \bar{\kappa}^r_1 = \frac{\bar{\kappa}_1
}{\|x\|}\|[A \,\,b]\|_F,
$$
see Theorems \ref{improveTLSthm} and \ref{myfirstestimate},
respectively. We calculate the approximate condition number used in
\cite{BjorckHeggernesMatstoms:2000}:
$$\bar{\kappa}^r_{\cite{BjorckHeggernesMatstoms:2000}}=\frac{\hat{\sigma}_1}{\hat{\sigma}_n -
\sigma_{n+1}}.
$$

As indicated by the table, all the given TLS problems are similar in
that $\sigma_{n+1}$ and $\hat{\sigma}_n$ are close but
$\sigma_{n+1}$ and $\sigma_n$ are not so close. As estimates of
$\kappa^r_g(A,b)$, the lower bounds $\underline{\kappa}^r_{2}$ are
very accurate, and the upper bounds $\bar{\kappa}^r_{2}$ improve the
corresponding $\bar{\kappa}^r_{1} $ and
$\bar{\kappa}^r_{\cite{BjorckHeggernesMatstoms:2000}} $
significantly by one or two orders.

{\bf{Example 2.}} In this example, the TLS problems are generated by
the function described in Appendix. For given $m, n$ and $\alpha$,
$A$ and $b$ are generated by
$$
[A \,\,b]= {\text{generate}}Ab\alpha(m,n, \alpha).
$$
A different $\alpha$ gives rise to a different TLS problem with
different properties. As $\alpha$ becomes smaller, $\sigma_{n+1}$
and $\hat{\sigma}_n$ become closer, so that the TLS problem becomes
worse conditioned. For each of the TLS problems, we calculate the
same quantities as those in Example 1 and list them in Table 2 in
which the first set of data is for $(m,n)=(500,350)$ and the second
set is for  $(m,n)=(1000,750)$.

\begin{table}[h]\label{tab2}
\begin{center}\caption{}
{ \doublerulesep18.0pt \tabcolsep0.01in
\begin{tabular}{|c |c| c | c| c | c |c|c |}\hline
$\alpha$   & $\sigma_{n+1}/\sigma_n$ &
$\sigma_{n+1}/\hat{\sigma}_n$& $\kappa^r_g(A,b)$ &
 $\underline{\kappa}^r_{2}$& $\bar{\kappa}^r_{2} $  & $\bar{\kappa}^r_{1} $ &
$\bar{\kappa}^r_{\cite{BjorckHeggernesMatstoms:2000}}$
 \\
\hline $10^{-2}$   & $0.953$ &
 $1-3.05\times 10^{-4}$ & $2.55\times 10^4$ & $8.98\times 10^3$&
 $1.60\times 10^5$&$5.14\times 10^5$ & $6.29\times 10^5$
 \\
$10^{-3}$   & $0.980$ &
 $1-3.16\times 10^{-6}$ & $2.01\times 10^5~$ & $8.75\times 10^4$&
 $2.42\times 10^6$&$4.92\times 10^7$ & $6.03\times 10^7$
 \\
$10^{-5}$ & $0.953$ &
 $1-2.77\times 10^{-10}$ & $1.97\times 10^7$ & $9.78\times 10^6$&
 $1.74\times 10^8$&$5.87\times 10^{11}$ & $ 7.20\times 10^{11}$
 \\
$10^{-7}$   & $0.966$ &
 $1-1.80\times 10^{-14}$ & $3.28\times 10^9$ & $1.12\times 10^9$&
 $2.38\times 10^{10}$&$8.38\times 10^{15}$ & $1.02\times 10^{16}$
 \\
  \hline
$10^{-2}$   & $0.983$ &
 $1- 2.78\times 10^{-4}$ & $6.76\times 10^4$ & $1.65\times 10^4$&
 $4.97\times 10^5$&$9.90\times 10^5$ & $1.21\times 10^6$
 \\
$10^{-3}$   & $0.978$ &
 $1-1.95\times 10^{-6}$ & $6.70\times 10^5$   & $1.93\times 10^5$&
 $5.09\times 10^6$&$1.38\times 10^8$ & $1.69\times 10^8$
 \\
$10^{-5}$   & $0.968$ &
 $1-3.01\times 10^{-10}$ & $4.33\times 10^7$   & $1.60\times 10^7$&
 $3.52\times 10^8$&$9.24\times 10^{11}$ & $1.13\times 10^{12}$
 \\
$10^{-7}$   & $0.993$ &
 $1-3.82\times 10^{-14}$ & $1.13\times 10^{10}$   & $1.44\times 10^9$
 & $7.02\times 10^{10}$&$7.38\times 10^{15}$ & $9.03\times 10^{15}$
 \\
 \hline
\end{tabular}}
\end{center}
\end{table}

We can see from the table that, for $\alpha = 10^{-2}$ in which
$\hat{\sigma}_n$ and $\sigma_{n+1}$ are not very close,
$\bar{\kappa}^r_1$ and
$\bar{\kappa}^r_{\cite{BjorckHeggernesMatstoms:2000}}$ are very
tight and they estimate $\kappa^r_g(A,b)$ quite accurately; for
$\alpha \leq 10^{-3}$, $\hat{\sigma}_n$ and $\sigma_{n+1}$ become
closer with decreasing $\alpha$, $\bar{\kappa}^r_1$ and
$\bar{\kappa}^r_{\cite{BjorckHeggernesMatstoms:2000}}$ estimate
$\kappa^r_g(A,b)$ increasingly more poorly. In contrast, however,
for all the cases, since ${\sigma}_n$ and $\sigma_{n+1}$ are not so
close, $\underline{\kappa}^r_{2}$ and $\bar{\kappa}^r_{2} $ estimate
$\kappa^r_g(A,b)$ accurately. Particularly, for $\alpha \leq
10^{-5}$, $\bar{\kappa}^r_{2} $ improves $\bar{\kappa}^r_1$ and
$\bar{\kappa}^r_{\cite{BjorckHeggernesMatstoms:2000}}$ very
considerably by several orders.

\section{Concluding Remarks}

In the paper, we have mainly studied the condition number of the TLS
problem and its lower and upper bounds that can be numerically
computed cheaply. For the TLS condition number, we have derived a new
closed formula. For a computational purpose, we can use it to
compute the condition number
more accurately. We have derived a few bounds, which are quite
sharp and can be calculated cheaply. We have confirmed our results
numerically and demonstrated the tightness of our bounds by
numerical experiments.

\bigskip
\bigskip
\leftline {\bf ACKNOWLEDGEMENTS} The work was partially supported by
National Basic Research Program of China 2011CB302400 and National
Science Foundation of China (No. 11071140) and Specialized Research
Fund for the Doctoral Program of Higher Education (No. 20070200009)

\appendix
\section{ Codes for generating tested TLS problems}

The following codes produce an $m \times (n+1)$ matrix $[A \,\,b]$,
which has the SVD $[A \,\,b]=U \Sigma V^T$ with $V(n+1,
n+1)=-\alpha$, where $0< \alpha <1$.

\begin{eqnarray}\nonumber
&& [A \,\,b] = {\text{{\bf{generate}}}}{\bf{Ab}}{\bf{\alpha}}(m,n,  \alpha) \\
\nonumber && \%~ m, n: {\text{two given positive integers with }} m \geq n \\
\nonumber
        && \%~ \alpha: {\text{a given positive number with $0<
\alpha <1$}} \\ \nonumber && \text{Generate}~ \tilde{V};~~ \text{\%
a random orthogonal matrix }~ \text{of order}~ n \\\nonumber &&V =
\text{\bf{generateV}} (n, \tilde{V}, \alpha);\\\nonumber &&B =
\text{rand}(m, n+1);~~ \text{\% the  Matlab function
rand(~)}\\\nonumber &&[U,\Sigma,\hat{V}]=\text{svd}(B, 0); ~~
\text{\% the  Matlab function svd(~)}\\\nonumber &&[A \,\,b] = U *
\Sigma
* V^T
\end{eqnarray}

The subfunction ${\bf{generateV}}( )$ is shown as follows. It is
used to produce an $(n+1) \times (n+1)$ orthogonal matrix $V$ with
 $V(n+1, n+1) = -\alpha$, where
$0< \alpha <1$. The idea of construction comes from Lemma
\ref{CSdecomp}.
\begin{eqnarray}\nonumber
    &&[V]={\text{{\bf{generate}}}}{\bf{V}}(n, \tilde{V}, \alpha) \\ \nonumber
    && \%~ n: {\text{a given positive integer}} \\ \nonumber
    && \%~ \tilde{V}: {\text{a given orthogonal matrix of order $n$}} \\ \nonumber
    && \%~ \alpha: {\text{a given positive number with $0<
\alpha <1$}} \\ \nonumber
  && \text{partition}~ \tilde{V} = [\tilde{v}_1, \ldots, \tilde{v}_n]; \\ \nonumber
  &&\text{generate}~ U=[u_1, \ldots, u_n];~~~  {\text{\% a random orthogonal
  matrix of order $n$}}\\ \nonumber
    && V_{11} = [u_1, \ldots, u_{n-1}] [\tilde{v}_1, \ldots, \tilde{v}_{n-1}]^T
    + \alpha u_n \tilde{v}^T_n; \\ \nonumber
  && V = \left[
           \begin{array}{cc}
             V_{11} & \sqrt{1-\alpha^2} u_n \\
             \sqrt{1-\alpha^2} \tilde{v}^T_n & -\alpha \\
           \end{array}
         \right]
\end{eqnarray}

{\bf Remark.} Lemma 4.3 in \cite{GolubVanLoan:1980} gives
\begin{equation}\nonumber
\frac{|\hat{u}^T_n b |}{2 (\hat{\sigma}_n - \sigma_{n+1})} \leq
\|x\| \leq \frac{ \|b\|}{\hat{\sigma}_n - \sigma_{n+1}}.
\end{equation}Equivalently, it holds that
\begin{equation}\label{Golubineqn}
\frac{|\hat{u}^T_n b |}{2 \|x\|} \leq \hat{\sigma}_n - \sigma_{n+1}
\leq \frac{ \|b\|}{ \|x\|},
\end{equation}where it is supposed that $x \neq 0$.
Note that $V(n+1, n+1)=-\alpha$ and $\alpha = \frac{1}{\sqrt{1+
\|x\|^2}}$. From (\ref{Golubineqn}) we see that a small $\alpha$
implies that $\hat{\sigma}_n$ and $\sigma_{n+1}$ are close in some
sense.

\newcommand{\Gathen}{\relax}\newcommand{\Hoeij}{\relax}

\end{document}